
\magnification\magstep1
\scrollmode
\hsize=13 cm
\vsize= 20cm
\overfullrule=0 pt

\font\tensmc=cmcsc10              

\font\tenbbb=msbm10 \font\sevenbbb=msbm7 \font\fivebbb=msbm5
\newfam\bbbfam
\textfont\bbbfam=\tenbbb \scriptfont\bbbfam=\sevenbbb
  \scriptscriptfont\bbbfam=\fivebbb
\def\Bbb{\fam\bbbfam}      
 
\font\tengoth=eufm10 \font\sevengoth=eufm7 \font\fivegoth=eufm5
\newfam\gothfam
\textfont\gothfam=\tengoth \scriptfont\gothfam=\sevengoth
\scriptscriptfont\gothfam=\fivegoth
\def\goth{\fam\gothfam}    
 
\def\opname#1{\mathop{\rm#1}\nolimits} 
 
\def\a{\alpha}                    
\def\A{{\cal A}}                  
\def \AU{{\cal A(U)}}             
\def\B{{\cal B}}                  
\def\bbuildrel#1_#2^#3{\mathrel{
     \mathop{\kern 1pt#1}\limits_{#2}^{#3}}}
\def\bw{\bigwedge\nolimits}       
\def\cinfty{C^\infty}             
\def\cite#1{\lbrack{\bf#1}\rbrack} 
\def\De{\Delta}                   
\def\dim{\opname{dim}}            
\def\ENE{{\cal N}}                
\def\eq#1{{\rm(#1)}}              
\def\FL{{\cal FL}}                
\def\frac#1#2{{#1\over  #2}}      
\def\Ga{\Gamma}                   
\def\harr#1#2{\smash{\mathop{\hbox to .5in{\rightarrowfill}}
     \limits^{\scriptstyle#1}_{\scriptstyle#2}}}
\def\hharr#1#2{\smash{\mathop{\hbox to .3in{\rightarrowfill}}
     \limits^{\scriptstyle#1}_{\scriptstyle#2}}}
\def\harrt#1#2{\smash{\mathop{\overline{\hbox to .5in{\hrulefill}}}
     \limits^{\scriptstyle#1}_{\scriptstyle#2}}}
\def\ID{\relax{1\kern-.24 em\rm l}}  
\def\id{\opname{id}}              
\def\Id{{\cal I}}                 
\def\Im{\opname{Im}}              
\def\L{{\cal L}}                  
\def\ker{\opname{ker}}            
\def\M{{\cal M}}                  
\def\o{\omega}                    
\def\Om{\Omega}                   
\def\ox{\otimes}                  
\def\o+{\oplus}                   
\def\O+{\bigoplus}                
\def\pd#1#2{{{\partial#1}\over{\partial#2}}}  
\def\R{{\Bbb R}}                  
\def\S{{\Bbb S}}                  
\def\set#1{\{\,#1\,\}}            
\def\smc{\tensmc}                 
    
\def\srule#1#2{(-1)^{\vert#1\vert\,\vert#2\vert}} 

\def\th{\theta}                   
\def\Th{\Theta}                   
\def\T{{\bf T}}                   
\def\Tau{{\cal T}}                
\def\Tk{{{\bf T}^{(k)}}}          
\def\Tkmu{{{\bf T}^{(k-1)}}}      
\def\Tu{{\bf T^{(1)}}}            
\def\U{{\cal U}}                  
\def\v{\vee}                      


\def\X{{\goth X}}                 
\def\Xk{{X^{(k)}}}                
\def\Xkmu{{X^{(k-1)}}}            
\def\Xl{{X^{(l)}}}                
\def\Yk{{Y^{(k)}}}                
\def\z{\zeta}                     
\def\1{\'{\i}}                    
\def\3{\sharp}                    
\def\7{\dagger}                   
\def\.{\cdot}                     
\def\:{\colon}                    
\def\<#1>{\langle#1\rangle}       


\def\today{\number\day\space      
	   \ifcase\month\or
	    January\or February\or
	    March\or April\or May\or
	    June\or July\or August\or
	    September\or October\or
	    November\or December\fi
	    \number\year}
	    
 \newcount\secnum
\outer\def\beginsection#1. #2\par{
      \vskip 0pt plus.1\vsize        
      \penalty -250                  
      \vskip 0pt plus-.1\vsize       
      \bigskip \vskip\parskip        
      \global\secnum=#1%
      \global\equationnum=0
      \message{#1. #2}%
      \leftline{\bf#1. #2}\nobreak
      \smallskip\noindent}
 
\outer\def\subsection#1.#2. #3\par{%
      \ifnum#2=1 \smallskip       
       \else \bigskip \fi         
      \message{Sec #1.#2.}%
      \leftline{\it#1.#2. #3}\nobreak
      \smallskip\noindent}
 
\def\declare#1. #2\par{\medskip   
	  \noindent{\bf#1.}\rm    
	  \enspace\ignorespaces
	  #2\par\smallskip}
 

\newcount\equationnum             
\def\ecnum{\the\secnum.
	   \the\equationnum}      
\everydisplay={\global\advance    
	       \equationnum by 1  
	       \leftnumdisplay}   
\def\leftnumdisplay#1{\let\nxxt#1\maybealign}
\def\maybealign{\ifx\nxxt\leqalignno
		\else\leftnumeqn\fi\nxxt}
\def\leftnumeqn#1$${#1\leqno(\ecnum)$$}  
 
\def\eqlabel#1{\xdef#1{\ecnum}}   
 
\def\refno#1. #2\par{\smallskip   
	  \item{\lbrack#1\rbrack}
		#2\par}
 
\long\def\suspend#1\resume{}      


\def\map#1#2#3{#1:#2\to #3}
\def\vectorfields#1{{\goth X} (#1)}
\def\forms#1#2{\bigwedge^#1(#2)}
\def\rank{\opname{rank}}       
 
\def\Marmuk{5}
\def\Grapon{4}
\def\MCS{6}    
\def\Pide{7}    
\def\Conv{1}
\def\Tulcz{8}    
\def\Tulczdos{9}       
\def\AAM{2}  

\def\Super{3}


\centerline{\bf SECTIONS ALONG MAPS IN GEOMETRY AND PHYSICS.}
 
\bigskip
 
\centerline{\smc Jos\'e F. Cari\~nena}

\smallskip
 
\centerline{\it Departamento de F\'{\i}sica Te\'orica,
Universidad de Zaragoza,}

\centerline{ \it 50009 Zaragoza, Spain.}
  
\bigskip
 
\vfill\vfill
Talk given at the meeting {\sl Geometry and Physics\/} on  the occasion
of the 65th birthday of  W.M. Tulczyjew 
 \eject
 
\beginsection 1. Introduction.

Geometric techniques have been applied to physics in many different
ways and they have provided powerful methods of dealing with classical
problems from a new geometric perspective. Vector fields, forms,
exterior differential calculus, Lie groups, fibre bundles, connections,
etc...,
are now well established tools in modern physics. One of the main
contributions of Prof. W.M. Tulczyjew is the alternative way of using
the geometric concepts  in a more algebraic approach, using, for instance,
derivations of algebras and related concepts. The  aim
of this paper is to point out some of my contributions 
 in this direction [\Conv--\Super,\MCS],
which have received a big influence of  Tulczyjew's works
[\Pide,\Tulcz,\Tulczdos]. Moreover, it
will be shown how this approach allows us to translate the usual concepts arising in
Geometrical Mechanics to the framework of Supermechanics.

\beginsection  2. Hamiltonian dynamical systems.

The geometric framework for the description of classical 
(and even quantum) systems 
is    the theory of Hamiltonian dynamical systems. They are triplets
$(M,\omega,H)$ where $M$ is a differentiable manifold,
$\omega\in Z^2(M)$ is a regular closed 2-form in $M$ and 
$H\in C^\infty (M)$ is a function called Hamiltonian. The dynamical vector field
$X_H$ is then the solution of the equation $i(X_H)\omega=dH$.

 We recall that an infinitesimal  symmetry of a  Hamiltonian dynamical
systems
 $(M,\omega,H)$ is  given by a Hamiltonian vector field $X\in \X(M)$ (i.e.,
 $i(X)\omega\in B^1(M)$) such that
 $XH=0$. Noether's theorem establishes a one-to-one correspondence between
 infinitesimal symmetries and constants of motion: 
 $XH=0$ and $i(X)\omega=df$\ if and only if  $f$ is a constant of
motion.

Very interesting  examples of HDS are those defined by regular
 Lagrangians, $(TQ, \omega_L, E_L)$,
with $\omega_L=-d\theta_L=-d(dL\circ S)$, $E_L=\Delta L-L$.
More accurately, the geometric approach to the Lagrangian description
makes use of the geometry of the tangent bundle of the
configuration space.

\subsection 2.1. The geometry of tangent bundles.

The  tangent bundle $\tau_Q\: TQ\to Q$ is characterized by
the existence of a vector field generating  dilations  along the fibres,
called Liouville vector field, 
  ${\Delta}  {\in }{ \goth X}(TQ)$,
  and 
the vertical endomorphism which is a $ (1,1)$--tensor field $S$ in $TQ$
 defined by $ S_{(q,v)}U={\xi }^{(q,v)}({\tau}_{*(q,v)}U), \,
\forall U\in T_{(q,v)}(TQ), 
 v\in T_qQ$,
 where
 ${\xi }^{(q,v)}\colon T_qQ\to T_{(q,v)}(TQ)$ denotes  the vertical lift
 defined by  $ 
 {\xi }^{(q,v)}(w) f= d/dt\left[ f(v+tw)\right]_ {\,\vert\,  t=0}$.  
In a natural coordinate system for $TQ$, induced from a chart in $Q$,
$ \Delta =v^i\partial/\partial v^i$,
and  $ S= (\partial/\partial v^i) \otimes dq{^i}.$

 The image under $S$ of a section for ${\tau}_ {TQ}\:T(TQ){\to} TQ$,
a vector field in $TQ$, is a new vector field in $TQ$,
 again.  This correspondence will also
be denoted $S$.

There is another vector bundle structure on $T(TQ)$ given by
$T\tau_Q:T(TQ){\to} TQ$. Vector fields on $TQ$ that are also sections for
$T\tau_Q$ are called SODE (second order differential equations).
They are characterized by $S(X)=\Delta$ and in tangent bundle coordinates
look like $ X= v{^i} (\partial/\partial q^i)+f{^i} (\partial/\partial v^i).$

A map $\varphi\:Q\to Q$ induces a map $\Phi=T\varphi\:TQ\to TQ$
such that $[T\Phi,S]=0$
and
$(T\Phi)(\Delta)=\Delta,$ and conversely,
if $\Phi\:TQ\to TQ$ satisfies these two properties,
then there exists a map $\varphi\:Q\to Q$ such that $\Phi=T\varphi$.
These transformations of the velocity
phase space are called point transformations.
Another characterization of these point tranformations is that if $X$ is a SODE, then
$T\Phi (X)$ is also a SODE. Notice that if $X$ is a SODE, then
$S(T\Phi(X))=T\Phi(S(X))=T\Phi(\Delta)=\Delta.$

\subsection 2.2. The Lagrangian formalism.

Given a function $L\in C^\infty (TQ)$, 
 we  define the 1--form  $\theta _L\in \bigwedge ^1(TQ)$ by
$\theta _L=dL\circ S$. When the exact 2--form $\omega _L= -d\theta _L$
 is nondegenerate the
 Lagrangain $L$ is called regular and then $(TQ,\omega_L)$ is a symplectic
 manifold. The energy function $E_L$ is given by $ E_L={\Delta  }(L) - L$.
The
 coordinate  expressions are
$ {\theta}_L= (\partial L/\partial v^i) dq{^i}$ and $ E_L=v^i(\partial L/\partial v^i)
 - L $.

 A remarkable fact is that if $\Phi$ is a point transformation, then
 $\Phi^*\omega_L=\omega_{\Phi^*L}$ and $ \Phi^*E_L=E_{\Phi^*L},$
 because
 $\Phi^*\theta_L=\phi^*(dL\circ S)=(dL\circ S)\circ (T\Phi)=
 (dL\circ S)\circ (T\Phi)=\Phi^*\theta_L$
 and
 $ \Phi^*E_L= \Phi^*(\Delta L- L)=\Delta (\Phi^*L)-\Phi^*L=E_{\Phi^*L}$.
However, for a general transformation $\Phi\:TQ\to TQ$, 
$\Phi^*\omega_L\not =\omega_{\Phi^*L}$ and $ \Phi^*E_L\not =E_{\Phi^*L}$.
 Therefore, the important point in the search of symmetries is that if $X=Y^c$ where
$Y\in \X(Q)$, then, $X$ is a symmetry of
$(TQ, \omega_L, E_L)$
iff $X$ is a symmetry of $L$ (up to a gauge term).
On the contrary, for  $X$ that are not   complete lifts,  symmetries of $L$
have nothing to do with  symmetries of the HDS $(TQ, \omega_L, E_L)$.
Complete lifts correspond to {\sl point transformations}.

In order to 
establish a one-to-one correspondence between symmetries of $L$
and constants of motion we shall generalize the concept of symmetry
in order to include ``non-point transformations". In the physicist's
language point transformations are written
$\delta q^i= \epsilon \,f^i(q)$ and 
$\delta v^i= \epsilon\, 
v^j\,(\partial f^i/\partial q^j),$
corresponding to the flow of $X=Y^c$ with 
$Y= f^i(q)(\partial/\partial {q^i})\in \X(Q)$.
Non-point transformations cannot be completed, 
$\delta q^i= \epsilon\, f^i(q,v)$, $ \, \delta v^i= ?$,
and this fact leads 
 to consider  ``objects" like    $Y=f^i(q,v)(\partial/\partial {q^i})$,
which are NOT 
 vector fields but vector fields along the tangent bundle 
projection, a concept that we will introduce in next section.
   
\subsection 2.3. Newtonoid vector fields.

 Marmo and Mukunda characterization of symmetries is as follows [\Marmuk]: 
Let $D$ be a SODE and ${\goth X}_D$ denote
$\X_D=\set{X\in\X(TQ)\mid S([X,D])=0}$.
There is a projection $\pi_{{}_D}\colon\X(TQ)\to
{\goth X}_D$, given by $\pi_{{}_D}(X)=X(D)=X+S([D,X])$.
In coordinates, if $X$
is  written $X=\eta^i\pd{}{q^i}+\xi^i\pd{}{v^i}$, then 
 $$X(D)=\eta ^i \pd{}{q^i}+(D \eta ^i)\pd{}{v^i}.$$

{\bf Theorem 1.} {\it 
Let $L\in\cinfty (TQ)$ be a regular Lagrangian. If $X\in\X(TQ)$ is
such that  $\exists F\in C^\infty(TM)$  satisfying}
$$      \L_{X(D)}L=\L_DF  \quad{\rm{{\it{for\ any}}\quad SODE\ }}D,
 \eqlabel\theor
$$
{\it then $G=i_X\theta_L-F$ is a constant of motion. Moreover, if $\Gamma$
is the 
dynamical vector field, then  ${\cal L}_{X(\Gamma)}\Omega_L=0$ and
 ${\cal L}_{X(\Gamma)}E_L=0$, i.e.
$X(\Gamma)$ is a symmetry of  $(TQ,\omega_L,E_L)$.

Conversely, if $X$ is a symmetry of $ (TQ,\omega_L,E_L)$, then $X=X(\Gamma)$
and there
exists a function $F\in \cinfty(TQ)$ such that~\eq\theor\ holds.}

When $X$ is a complete lift $ X=Y^c$ with 
$Y\in\X (Q)$,  i.e., the flow of $X$ is the differential of the flow of $Y$,
 then $X(D)=X$  for any SODE  $D$, and $F$ reduces to the
pullback $\tau^\star h$ of a function $h$ on the base $Q$.

If a  vector field $X$ of $\X(\Gamma)$  is a
symmetry of $(TQ,\omega_L,E_L)$,  then its vertical components
 $\xi^i$ are
determined by the other ones:  $\xi^i=\Gamma \eta^i$. Essentially we should consider
{\sl equivalence classes} of  vector fields
$$[X]=X+{\rm Vert}.$$
\beginsection 3. Sections along maps.

Let  $\map{\pi  }{E}{M}$ be a fibre bundle and $\map{\phi   }{N}{M}$ a
differentiable map. A section along $\phi $ is a map $\map{\sigma}{N}{E}$
such that $\pi \circ   \sigma  =\phi   $.
They are in a 1-1 correspondence with sections of the induced bundle $\phi  
^*E$,
$$
\phi   ^*E=\set{(n,e)\in N\times E\vert\,\phi
   (n)=\pi  (e)}\subset
N\times E.
$$

The set of sections along $\phi$ will be noted 
$\Sigma_{\phi}(\pi)$. When
 $E$
is a vector bundle the set $\Sigma_{\phi}(\pi)$ is endowed
 with a
 $\cinfty(N)$--module structure. For more details see [\Conv,\AAM]. 

In particular we will be interested in the
vector bundles
$\tau_M:TM\to M$, $\pi^p_M:(T^*M)^{\wedge p}\to M$ and $\rho^p_M:(T^*M)^{\wedge p}\otimes
TM\to M$,
and in these cases we will denote 
$\X(\phi)=\Sigma_{\phi}(\tau_M)$,
$\bigwedge\,^p
(\phi)=\Sigma_{\phi}(\pi _M^p)$ and $
V^{p}(\phi)=\Sigma_{\phi   }{\rho_M^p}$,
respectively. When $N=M$ and $\phi=\id$ the set $\X(\id)$ coincides with 
$\X(M)$ and the set $\forms{p}{\id   }$ reduces to 
$\forms{p}{M  }$.
 \medskip
{\bf Examples:}

Let  $\map{\gamma}{{\Bbb R}}{M}$ be a curve in $M$.  The tangent 
vectors $\dot \gamma$
define a section  $\map{\dot \gamma   }{{\Bbb R}}{TM}$ of $\tau_M$ 
along  $\gamma  $. The restriction  of $X\in \vectorfields M$
on the curve $\gamma$ is also a vector
field along $\gamma$.

The generalization of these examples is:
Let $\phi   $ be a map from $N$ to $M$. A vector field  $Y\in \vectorfields N$
defines a vector field  along $\phi$ by  
 $T \phi \circ Y\in \vectorfields{\phi   }$. Similarly, when  $X\in 
\vectorfields M$
the 
restriction $X\circ \phi   $ of $X$ on the  image by $\phi   $ is a vector
field along $\phi   $. The above vector fields  $X$ and 
$Y$ are said to be  $\phi   $-related when  $X\circ \phi   $ and
$T{\phi  }{}\circ Y$  coincide along $\phi   $. 
Similarly, if $\beta $ is a $p$-form in $M$, the restriction  $\beta \circ \phi  
$ of $\beta $ on the  image by  $\phi   $ is a $p$-form along 
$\phi   $.
\medskip

  Given  
$\alpha  \in \forms{p}{\phi   }$,  $T^*{\phi  
}{}\circ \alpha  $ is a $p$-form in $N$. The pull-back by $\phi$ of
$\beta \in\forms{p}{M}$ is obtained by iteration of both processes $\phi  
^*(\beta )=T^*{\phi   }{}\circ \beta  \circ \phi   $.
 
 When  $E$ is a vector bundle and $\{\sigma
_\alpha  \}$ is a local basis  of  $\Sigma(\pi)$,
 then  $\{\sigma _\alpha 
\circ \phi   \}$ is a local basis  of   $\Sigma_{\phi}(\pi)$,
and  $\sigma  \in \Sigma_{\pi  }$ can be written as  $
\sigma  =\zeta^\alpha  (\sigma  _\alpha  \circ \phi   )$ with
$\zeta^\alpha 
\in\cinfty(N)$.

In the above case, taking local  coordinates  
$(z^A)$ in $N$ and $(x^i)$ in $M$ we have
$$
\eqalign{
&X\in \vectorfields{\phi   }\qquad
 X=X^i\left(\pd{}{x^i}\circ \phi   \right)\cr
&\alpha  \in \forms{p}{\phi   }\qquad
 \alpha  =\alpha  _{i_1\ldots i_p}(dx^{i_1}\circ \phi   )\wedge   \ldots \wedge  
(dx^{i_p}\circ \phi   )} 
$$
where  $X^i$  and  $\alpha  _{i_1\ldots i_p}$ are functions in $N$.
When $N=TM$ and $\phi$ is the projection $\tau_M$ the vector fields and forms
along   $\tau_M$ are written 
$$X=X^i(x,v) \left(\pd{}{x^i}\circ\tau_M  \right),\quad 
\alpha  =\alpha  _{i_1\ldots i_p}(x,v) (d{x^{i_1}}\circ \tau_M   )\wedge   \ldots \wedge  
(d{x^{i_p}}\circ \tau_M   )
$$

Vector fields along $ \phi $  act on  functions on $M$ giving rise to
functions on $N$. 

If $X\in\vectorfields  \phi $ and $n\in N$ then $X(n)$ is a tangent vector to
$M$ at the point $ \phi (n)$ which  acts on a function $h\in\cinfty (M)$ by 
$(Xh)(n)=X(n)h$. The Leibnitz rule for tangent
vectors implies that $$X(hl)=\phi^*h\,Xl+\phi^*l\,Xh.
$$

 A map satisfying this
property is called a $\phi^*$-derivation (of degree 0) [\Pide].

{\bf Definition 1.} Let $\phi:N\to M$
a differentiable map. A $\phi^*$-derivation of degree $r$ 
of scalar
  on $M$  is a $\Bbb R$-linear map $ D:\bigwedge(M)
\to \bigwedge(N)$ satisfying 
$$
D\left(\bigwedge^p(M)\right)\subset \bigwedge^{p+r}(N),  \quad
 D(\alpha\wedge\beta)=D\alpha\wedge \phi^*\beta+(-1)^{pr}
\phi^*\alpha\wedge D\beta
$$
for $\beta\in\forms qM$ and $\alpha\in\forms pM$. It is said to be of type 
$i_*$ when $D g=0, \, \forall g\in C^\infty (M)$.

For instance, given a vector field along $\phi:N\to M$, 
$X$,  a type $i_*$ $\phi$-derivation
$i_X:\bigwedge ^p(M)\to \bigwedge^{p-1}(N)$ of degree $-1$ is defined by 
$i_X g=0 \, \forall g\in C^\infty (M)$  and 
$$
(i_X \omega)_z(v_1,\dots,v_{p-1})=
\omega_{\phi(z)}(X_z,\phi_{*z}v_1,\ldots,\phi_{*z}v_{p-1})
$$
where $v_1,\ldots,v_{p-1}\in T_z(N)$.

By a type $d_*$ $\phi$-derivation of degree $r$ we mean that $D\circ  d_{(M)}
=(-1)^r
 d_{(N)}\circ D$. An example 
of such a type 
$\phi$-derivation, $d_X$, is defined by 
$$
d_X=i_X\circ d_{(M)}+d_{(N)}\circ i_X,
$$
where $d_{(M)}$ stands for the operator of exterior differentiation in $M$.
 This
is of type $d_*$, i.e., $d_X\circ d_{(M)}=d_{(N)}\circ d_X$.

Note that when $X\in\X(\id_M)\equiv\X(M)$ the $\id_M$-derivations $i_X$ and
$d_X$ are but the contraction or inner product $i(X)$ (or $i_X$) and the Lie
derivative ${\cal L}_X$, respectively. For this reason, $i_X$ and $d_X$ will be
called {\it contraction} and {\it Lie derivative}, respectively.

There exists a  section along $\pi$ in each vector bundle  $\map \pi 
EM$, which is given by the identity map in $E$.
When choosing local coordinates $(x^i,y^\alpha  )$ in  $E$  and a local
basis $\{\sigma  _\alpha  \}$ of  sections for $\pi:E\to M$ such that  $y^\alpha 
(e)=\sigma  _\alpha  (\pi  (e))$, for  $e\in E$, then the local
expression  of  ${\cal C}$ is
$ {\cal C}=y^\alpha  (\sigma  _\alpha  \circ \pi 
)$.

The most important cases in Classical Mechanics are those of $E=TM$ or $
E=T^*M$. Then $\cal C$ reduces in these cases to the \lq \lq total time
derivative" $\bf T$ (in the time-independent formalism) and the Liouville
1-form, to be denoted $\check{\theta    }_0$, up to an identification of
$\pi$-semibasic forms with forms along $\pi$ (If  $\phi   $ is a submersion,
every $p$-form along $\phi   $ may be identified to a $\phi$-semibasic
$p$-form in $N$.)

The  coordinate expressions for ${\bf T}$ and $\check{\theta}_0$ are
$$
{\T}= v^i\left(\pd{}{x^i}\circ \tau _M\right)
\qquad{\rm {and}}\qquad
\check{\theta }_0=p_i(dx^i\circ \pi_M).
$$

\beginsection 4. Applications in Geometry.

\subsection 4.1. The geometry of $TQ$ revisited.

The fundamental objects in tangent bundles can be introduced
in an alternative way that can be
generalized for the case of Supergeometry and Supermechanics.

If $f\in C^\infty (Q)$, let $f^V\in \cinfty (TQ)$ be defined by
$$f^V=d_{\bf T} f:= \sum^m_{i=1} \pd{F}{q^i}v^i ,
$$
where $F:= \tau^*(f)$.

A  vector field $Y$ on $TQ$ is determined
by its action on the functions $f^V$: 
if $Y\in\X(TQ)$ satisfies $Y(f^V) = 0,  \, \forall f\in\X(TQ)$,
then $Y\equiv 0$.
This property allows us to define the vertical lift:
If $X\in\X(Q)$, then $X^V\in\X(TQ)$ is defined by:
$X^V(f^V) = \tau^*\bigr(X(f)\bigl), 
\, \forall f\in\cinfty (Q)$.

Similarly, we can define the vertical lift $X^V\in\X(TQ)$
of a vector field along $\tau$, $X$, 
 by the relations $
X^V(f^V) = X(f), 
\, \forall f\in\cinfty (Q)$.

The Liouville vector field $\Delta\in \X(TQ)$ can then be defined as the vertical
lift of the total time derivative $\bf T$: $\De= \T^V$. 
 On the other hand, the  vertical endomorphism is the
$(1,1)$--tensor field
$S \: \X(TQ) \to \X(TQ)$  defined by $S(Y) := T\tau(Y)^V$.

 \subsection 4.2. Covariant derivative of a vector field along a curve.

 A curve $\gamma\:\R\to M$ has associated a vector field along $\gamma$,
 the tangent vector field $\dot \gamma$.
In particular, if $X\:\R\to TQ$ is a
 vector field along the curve $\sigma\:\R\to Q$, it can be considered as
a curve in $TQ$.  The associated vector field along $X$, $\dot X$, can be
composed with the natural isomorphism $\Psi\:T(TQ)\to T(TQ)$, that in local
coordinates is given by $\Psi(q,v,\dot q,\dot v)=(q,\dot q,v,\dot v)$,
giving rise to a vector field along $\dot \sigma$, $X^1=\Psi\circ \dot X$.

Given a connection on $Q$, we can also lift horizontally $X$ and we obtain
the vector field along $\dot\sigma$, $X^H$, given by $X^H(t)=\xi^H_{\dot \sigma(t)
}(X(t))$.
The difference $X^1-X^H$ is a vertical vector field along $\dot \sigma$,
and therefore, there exists a vector field along $\sigma$,
to be denoted $DX/Dt$,
such that
$$
(X^1-X^H)(t)=\xi^V_{\dot \sigma(t)}\left(\frac{DX}{Dt}\right).
$$

The vector field along $\sigma$ $DX/Dt$ is called total covariant
 derivative of $X\in \X(\sigma)$.
So, if
$X(t)=\eta^i(t)(\partial/\partial {x^i})_{\vert  \sigma(t)},$
 then
$$
 X^1(t)=\eta^i(t)(\pd{} {x^i})_{\vert \dot \sigma(t)}
+\frac{d\eta^i}{dt}(t)
(\pd{} {v^i})_{\vert \dot \sigma(t)},$$
 and
 $$ X^H(t)=\eta^i(t)(\pd{} {x^i})_{\vert \dot \sigma(t)}-
\Gamma^i_j(\dot\sigma(t))  \eta^j(t)(\pd{} {v^i})_{\vert \dot \sigma(t)}.
 $$
Therefore,
 $$
 \frac{DX}{Dt}(t)=\left[\frac{d\eta^i}{dt}(t)+\Gamma^i_j(\dot\sigma(t))
 \eta^j(t)\right]\pd{}{x^i}_{\vert \dot \sigma(t)}.
$$

This allows us to define the concepts of parallelism and geodesic
curves in the well known way. The equation for geodesic curves is:
$$\frac{d^2\sigma^i}{dt^2}+\Gamma^i_j(\dot\sigma(t))\frac {d\sigma^j}{dt}=0.
$$

\beginsection 5. Applications in Physics.

\subsection 5.1. Generalized symmetries.

  If $k,l\in {\Bbb N}$, there is a natural  inmersion
$i_{k,l}:T^{k+l}Q\to T^l(T^kQ),$
given by
$[\rho]^{k+l}\mapsto [\tilde \rho^k]^l$.
 We remark that  for $k=0, l=1$ we obtain the identity in TQ, i.e.,
$i_{0,1}={\bf T}$. In general $i_{k,1}:T^{k+1}\to T(T^kQ)$ is such that
$\tau_{T^kQ}\circ i_{k,1}=\tau_{k+1,k}$,
and therefore $i_{k,1}$ is a vector field along  $\tau_{k+1,k}$
which will be denoted ${\bf T}^{(k)}$.

If  $X$ is a vector field along $\tau_{1,0}:TQ\to Q$, then 
there exists one
vector field along $\tau_{2,1}:T^2Q\to TQ$, denoted $X^{(1)}$,  such that
$X \circ\tau_{2,1}= \tau_{1,0*}\circ X^{(1)}     $
and satisfying the  commutation property  $d_{ X^{(1)} }\circ d_{T^{(0)}}=
 d_{ T^{(1)} }\circ d_{X}$.

Second order differential equations $\ddot q=F(q,\dot q)$ can be seen
not only as vector fields $\Gamma\in \X(TQ)$, but 
as sections $\gamma$ of $\tau_{2,1}$, i.e., $\tau_{2,1}\circ\gamma=
\id_{TQ}$,  given by $a=F(q,v)$.

The relation between the two alternatives is given trough the
{\sl time derivative} vector field along $\tau_{2,1}$:
$\Gamma={\bf T}^{(1)}\circ \gamma$,
or in other words, ${\cal L}_{\Gamma}=\gamma^*\circ d_{{\bf T}^{(1)}}$.

In a similar way there exists a one-to-one correspondence
$I_\Gamma:\X(\tau_{1,0})\to \X(\Gamma)$, given by
$X\mapsto X^{(1)}\circ \gamma$,
which is but the  inverse of the restriction of $\tau_{1,0*}$ onto
$\X(\Gamma)$.

Moreover, it is then possible to show the following theorem [\Conv]:

{\bf Theorem 2.} {\it  Let $L$ be a regular Lagrangian, $\Gamma$ the
dynamical vector field satisfying
$i(\Gamma)\omega_L=dE_L$,
and $\gamma:TQ\to T^2Q$ the corresponding section.
If $X$ is  a vector field along $\tau_{1,0}$ 
and   there exists a function $F\in C^\infty(TM)$ such
that 
$  X^{(1)}  L= d_{ T^{(1)} }F$,
then the function $G=F-\check\theta_L(X)$ is a constant of motion.
The vector field $X^{(1)}\circ \gamma=X(\Gamma)$ is a symmetry of $\Gamma$.

 Conversely
if $G$ is a first integral of the motion given by $L$, then there exist 
$X\in\X(\tau_{1,0})$ and $F\in C^\infty(TM)$ such that the above
relation holds.}

Therefore if we call symmetries of $L$ to those vector fields 
$X\in \X(\tau_{1,0})$ satisfying the preceding condition, then there will 
be a one-to-one  correspondence between {\sl generalized infinitesimal 
symmetries} of $L$ and constants of motion. 
Point symmetries correspond to $X=Y\circ\tau_{1,0}$, with $Y\in \X(Q)$ 
and then $F$ is a basic function.

\subsection 5.2. The evolution $K_L$-operator.

The Legendre map $\FL:TM \to T^*M$, as well as the  time evolution operator
$K_L:\cinfty(T^*M)\to\cinfty(TM)$,  can also be defined as sections
along maps [\Grapon,\AAM]:
The 1-form $\theta _L$ is semibasic and can be
seen as a 1-form along $\tau$, $\theta_L\in\bigwedge(\tau)$, $\theta_L =
\pd{L}{v^i} (dq^i\circ \tau )$.
It is identified in this way with
the Legendre map.
Let $\chi $ be the natural diffeomorphism between $T^*(TM)$ and $T(T^*M)$,
with coordinate expresion  $\chi (x,v,p_x,p_v)=(x,p_v,v,p_x)$.
Then 
$K=\chi \circ  dL$ maps $TM$ in 
 $T(T^*M)$ in such a way that 
$\tau   _{T^*M}\circ  K_L={\cal F}L$, say, $K_L$ is a vector field along 
${\cal F}L$.

   In coordinates, $K_L$  is the  vector field along $ {\cal F}L$ given by
$$
K_L=v^i\left(\pd{}{x^i}\circ  {\cal F}L\right)
  +\pd{L}{x^i}\left(\pd{}{p_i}\circ  {\cal F}L\right)
$$
and it is very useful to relate constraint functions arising in the Hamiltonian
and Lagrangian formulations respectively.
It is determined by the following two equations:
$$i_{K_L}\omega _0=dE_L,\qquad T\pi_Q\circ K=\id_{TQ}.
$$

The main properties are that 
when applied to a constraint function in the Hamiltonian formalism,
it produces a constraint in the Lagrangian formalism, either a
dynamical constraint or even a SODE constraint. More especifically, 
when applied to a first class constraint function in the Hamiltonian
formalism produces a dynamical constraint that is ${\cal F}L$--projectable,
while when it is applied to a second class constraint function in the Hamiltonian
formalism produces a SODE constraint that is not ${\cal F}L$--projectable.

\subsection 5.3.  Applications in degenerate systems.

Let $(N,\omega )$ be a symplectic manifold and 
$\phi \:P\to N$ of constant rank. Given a  1-form $\alpha   $ in $P$, we
look for the set of points in which a solution of $i_\Gamma (\phi  ^*
\omega )=\alpha$
 exists, 
where  $\Gamma$ is a vector field in  $P$, i.e., we
are interested in the submanifold
$i_C\: C\to P$ of $P$ in which such a solution 
$\Gamma'\in     \X(C)$ does exist, namely $i_{\Gamma'}((\phi  \circ  i_C)^*\omega
)=i_C^*\alpha   $.

The above problem may be splitted in two. First we study the conditions for the
existence of  $X\in \X(\phi)  $ such that  $i_X\omega =\alpha$,
and then we determine the conditions for  $X$ to be image under $T\phi$ of
a vector field in  $P$. This is equivalent to the original problem because of
the relation
 $i_{T\phi\circ\Gamma}\omega =i_\Gamma(\phi  ^*\omega )$.
The second step is
but the condition for the solution to be tangent to 
$P$.
Using a well-known result of Linear Algebra we obtain that 
the equation  $i_X\omega =\alpha   $ has a solution with  $X\in   
\X(\phi)$ iff $p\in    P$  satisfies 
$$
\langle      z,\alpha   (p)\rangle      =0\qquad\quad{\rm{for all}} \
z\in    T_pP \quad {\rm such\ that\ }  T_{p}{\phi  }(z)=0.
$$

If 
 $X$ is a solution and  $Z\in    \X(\phi)$ is such that
$\hat{\omega }(\phi  (p))(Z(p))\in    \ker T^*_{p}(\phi)\quad
 \forall p\in    P$,
then  $X+Z$ is a solution too.

When 
$\alpha   $ is exact, $\alpha   =d F$, if for  $n\in   
\Im{\phi  }\subset N$ the submanifold  $\phi  ^{-1}(n)$ is connected,
then the above condition is equivalent to  $F$ to be  $\phi  $-projectable,
namely   there exists  $\widetilde F\in    \cinfty(N)$ such that
$\phi  ^*(\widetilde F)=F$.

The generalization to the case of a presymplectic manifold is: the equation 
$i_X\omega =\alpha   $ admits a solution iff   $\langle    
  z,\alpha   (p)\rangle      =0$ for any $z\in T_pP$ such that $T_{p}{\phi  }
  (z)\in 
   {\rm rad}(\omega )$, where $
{\rm rad}(\omega )=\set{v\in    TN\mid \omega (v,w)=0\quad \forall
  \omega \in    TN}$.

Once the condition holds in $P$, we look for 
the existence of a vector field  $\Gamma$
in $P$ such that 
$T{\phi  }\circ     \Gamma=X$,
which has a solution iff the equation   $T_{p}{\phi  }(\Gamma(p))=X(p)$ has
solution for any  $p\in    P$.

This is equivalent to 
 $\langle      X(p),\lambda   \rangle      =0, 
\,\forall \lambda  \in T^*_{\phi  (p)}N$ such that  $T^*_{p}{\phi  }(\lambda   )=0$,
or in other words, iff
$\langle      \delta    ,X\rangle      =0$,  for all $ \delta  
\in    \bigwedge^1(\phi)$ such that $T^*(\phi)\circ
\delta    =0$.
This gives rise to an inmersed  submanifold  $i_1 \: P_1\to P$ of $P$
and we repeat the preceding steps. 

If the image by  $\phi  $ is an inmersed submanifold  $j\: N_0\to N$ of $N$,
then a similar algorithm is used for finding a solution in  $N_0$.  If $\zeta
$ is a constraint function for  $N$, then  $\phi  ^*\zeta     $ is a  constraint
function for $P$. 

This is a generalization of what happens with the theory defined by a 
singular
Lagrangian when  $\phi  $ is the Legendre transformation.

 \subsection 5.4. Control systems. 

{\bf Definition 2.} Given a vector field along $\phi:N\to M$, we will 
say that a curve $\gamma:{\Bbb R}\to N$ is an integral curve 
for $X\in\X(\phi)$ if
$$X\circ \gamma=T\phi\circ  \gamma=(\phi\circ\gamma)^..
$$

If we choose coordinates $\{z^\alpha\}$ in $N$ and $\{x^i\}$ in $M$,
the vector field is
$X=X^i(z^\alpha)[(\partial/\partial {x^i})\circ \phi],$
and the integral curves are to be determined by solutions of the equation
$$\pd{\phi^i}{z^\alpha}(\gamma(t))\, \frac{d\gamma^\alpha}{dt}
=X^i((\gamma(t)).
$$
This system is not in normal form and the theorems of existence and uniqueness of solution
do not apply.

Let us consider the differential equation system
$$\frac{dx^i}{dt}=F(x^i,u^\alpha), \quad i=1,\ldots, n,
\quad \alpha =1, \ldots,m.
$$

>From the geometrical viewpoint this system can be seen as the
one determining
the integral curves of the vector field along the projection $\pi: M
\times \R^m\to M$,
where $\{x^i, i=1,\dots,x^n \}$ are the local coordinates of a point
in $M$,  
and  $\{u^\alpha, \, \alpha=1,\ldots ,m\}$
are the so called control functions. As a straightforward generalizations of this,
given a fibre bundle $\pi:B\to M$, a control system in $B$ is a vector field along
$\pi$, $X\in \X(\pi)$. A solution of the control system  is an
integral curve of  $X\in \X(\pi)$.

One of the main problems in the theory of control systems
 is to investigate
the set of points accesible  from one given point $p\in M$. More
especifically, a
control system is said to be controlable if for any given initial
point $p$
there exists an integral curve of the corresponding vector field
along $\pi$
such that $(\pi\circ \gamma)(0)=p$ and a value $t_1$ of the parameter
of the  curve $\gamma$ such that $(\pi\circ \gamma)(t_1)=q$.

The simplest example is when $M=\R^n$, $B=\R^n\times \R^m$ and the
equations describing the system are  of the linear type
$$\dot x^i=A^i\, _jx^j+B^i\,_\alpha u^\alpha,$$
where $A^i\, _j$ and  $B^i\,_\alpha$ are constant matrices. In this case
the Kalman rank controllability condition is that
$\rank (B, AB,\ldots, A^{n-1})=n.$

Another interesting example is when $X=u^1X_1(x)+\cdots+u^rX_r(x)$.
In this case, if the distribution ${\cal D}$ generated by  the vector
 fields
$X_1,\ldots,X_r$ is integrable, then the only accesible points from a
given point $x_0$ are those of its leaf. Otherwise, we should consider the
minimal integrable distribution containig  ${\cal D}$ and
then the system is controlable
if and only if $\dim D=n$.

\beginsection 6. Supermechanics.

The algebraic tools considered in the previous sections can be
translated to the framework of
Supergeometry, the theory of supermanifolds, and  Supermechanics [\Super].

 We recall that   a  graded manifold $\M$ is a pair $(M,\A)$, where
$M$ is a nice topological  space and $\A$ is a sheaf of superalgebras
over $M$ such that there are open sets $\U$ that cover $M$ such that
 $\AU\cong C^\infty(U) \ox \bw (\R^n)$ and satisfying
 glueing conditions on the overlaps.

Given $\Phi = (\phi,\phi^*)\: (\ENE,\B)\to (M,\A)$ a
supervector field along $\Phi$ is a morphism of sheaves over
$M$, $X\: \A \to \Phi_*\B$ such that 
$$
X(fg) = X(f)\, \phi^*_\U(g) 
+ \srule{X}{f} \phi^*_\U(f) \, X(g).
$$

As an example,  if $X\in \X(\A)$, then $\phi^*\circ X\in \X(\Phi)$, and 
 if $Y\in \X(\B)$ then $Y\circ\phi^*\in \X(\Phi)$.

The supertangent bundle $T^{k+1}\M$ is defined in a recursive way
 as follows. One has two morphisms: first,
 $\Tau_{k,k-1}\:T^k\M\to T^{k-1}\M$, and the corresponding
tangent map 
$T\Tau_{k,k-1}\:T(T^k\M)\to T(T^{k-1}\M)$; second, 
 $I_k\:T^k\M\to T(T^{k-1}\M)$ and
$\Tau_k\:T(T^k\M)\to T^k\M$ give rise to
$I_k\circ\Tau_k\:T(T^k\M)\to T(T^{k-1}\M)$.

Then  $T^{k+1}\M$ is the subsupermanifold of
$T(T^k\M)$ associated to the superideal
$$
\Id_{k+1}= <\set{\tau_k^*\circ i^*_k(F)-(T\tau_{k,k-1})^*(F)\:
F\in T(T^{k-1}\A)}>
$$

\noindent If $q^i_{(0,k)}$, $\dots,$ $q^i_{(k,k)},$ $v^i_{(0,k)},$
$\dots,v^i_{(k,k)},$ $\th^\a_{(0,k)},\dots,$
$\th^\a_{(k,k)},$ $\z^\a_{(0,k)},$ $\dots,$ $\z^\a_{(k,k)}$ are
supercoordinates on $T(T^k\M)$ then
$$
\Id_{k+1}= 
<q^i_{(1,k)}-v^i_{(0,k)},\dots,
q^i_{(k,k)}-v^i_{(k-1,k)},
\th^\a_{(1,k)}-\z^\a_{(0,k)},\dots,
\th^\a_{(k,k)}-\z^\a_{(k-1,k)}>
$$
and the supercoordinates in $T^{k+1}\M$ are $
\set{q^i_{0,k+1},\dots,q^i_{k+1,k+1},
\th^\a_{0,k+1},\dots,\th^\a_{k+1,k+1}}$.

The total derivative with respect to time
$\Tk\:T^k\A\to T^{k+1}\A$ is the element of 
$\X(\Tau_{k+1,k})$ such that $
\Tk(q^i_{j-1,k}) = q^i_{j,k+1}$, and
$\Tk(\th^\a_{j-1,k}) = \th^\a_{j,k+1}$,
for $j=1,\dots,k+1$. So, if $f\in \A$ then $
f^k_j := \tau_{k,j}^*\circ\T^{(j-1)}\circ\dots
\circ\Tu\circ\T(f)\in T^k\A$,
and if $Y\in\X(\A)$, the complete lift of $Y$ is
 $\Yk(f^k_j) =\bigl(Y(f)\bigr)^k_j,\quad(\forall 
f)(\forall j)$.

When  $X\in\X(\Tau_{k,0})$, the $l$--prolongation
 of $X$ is
$$
\Xl(f^l_j):= i^*_{k,l}\bigl((Xf)^l_j\bigr) \quad(\forall 
f)(\forall j).
$$

If $f\in T^{k-1}\A$ and $F=\tau^*_{k,k-1}f$
$$
f^V := 
\sum^m_{i=1}\sum^{k-1}_{j=0} 
{1\over j+1}\pd{F}{q^i_{j,k}}q^i_{j+1,k} 
+ \sum^n_{\a=1}\sum^{k-1}_{j=0}
{1\over j+1} \pd{F}{\th^\a_{j,k}}\th^\a_{j+1,k},
$$
is a superfunction in $T^k\A$.
If $Y\in\X(\A)$, the vertical lift of $Y$ is
$$
Y^V(f^V) = \tau^*_{k,k-1}\bigr(Y(f)\bigl) 
\qquad \forall f\in T^{k-1}\A.
$$

Moreover, if $X\in\X(\Tau_{k,k-1})$, the vertical lift
of $X$ is given by $X^V(f^V) = X(f), \, \forall f\in T^{k-1}\A.$
and the vertical superendomorphism is the
graded tensor of type $(1,1)$
$S_k \: \X(T^k\A) \to \X(T^k\A)$  defined by
$S_k(Y) :=(Y\circ\tau^*_{k,k-1})^V$.

On the other hand, he Liouville supervector field is $\De_k:=(\Tkmu)^V$, 
and the Cartan operator 
$\S^{(k)}\:\Om^1(T^k\A)\to \Om^1(T^{2k-1}\A)$  
is defined by
$$
\S^{(k)}:=\sum^k_{l=1}{(-1)^{l+1}\over l!}
\Tau^*_{2k-1,k+l-1}\circ d^{l-1}_{\T^{(k)}}\circ 
S^{*l}_k.
$$

A superdifferential equation of order $k+1$
is a $\Ga\in\X(T^k\A)$ such that $
\Ga \circ \tau^*_{k,k-1} = \Tkmu$, or in other words, 
$S_k(\Ga)=\De_k$.

There is a 1--1 correspondence between the
$\Ga$'s and the morphisms $\gamma\: T^{k+1}\AU \to T^k\AU$ 
such that $ 
\gamma \circ\tau^*_{k,k-1} = \id_{T^k\AU}$.

 The Cartan 1--form  defined by a super--Lagrangian $L\in T^k\A$
is
$\Theta_L:=\S^{(k)}(d\,L)$ and the Cartan 2--form is
$\Om_L:= -d\,\Th_L$.

If $k=1$,  $E_L:=\De L -L$, and then we say that
$\Ga\in\X(T^{2k-1}\A)$ is Lagrangian if
$\Ga$ is a superdifferential equation of order $2k$ such that 
$i_\Ga\Om_L= d\,E_L$.

We will say that
$X\in\X(\Tau_{2k-1,0})$ is a generalized
infinitesimal supersymmetry of the dynamical system
$(T\M,\Om_L,E_L)$ if there exists $F\in T^{3k-2}\A$
such that
$$
\Xk L = \T^{(3k-2)} F.
$$

The main theorem is the following:

 {\bf Theorem 3.} {\it Let  $L$ be a regular Lagrangian.
Then, if $X$ is a generalized supersymmetry there exists a
constant of motion $G$ such that
$$
\tau^*_{3k-1,2k-1}G= \<\Xkmu,
(\Tau^*_{3k-2,2k-1}\Th_L)^\v> -F, \eqlabel\symcon
$$
and conversely, if $G$ is a constant of motion,
there exists $X\in \X(\Tau_{2k-1,0})$ such that
$$
F:=\<\Xkmu,(\Tau^*_{3k-2,2k-1}\Th_L)^\v> - 
\tau^*_{3k-1,2k-1}G
$$
satisfies \eq\symcon}.

 \centerline{\bf References}
 {\narrower
\baselineskip=9.5pt

\medskip

[\Conv] {\smc Cari\~nena J.F., L\'opez C. and Mart\'{\i}nez E.},
{\sl A new approach to the converse of Noether's theorem},
{J. Phys. A: Math. Gen. {\bf 22} (1989), pp. 4777--87.}

 \medskip

[\AAM] {\smc Cari\~nena J.F., L\'opez C.  and  Mart\'{\i}nez E.},
{\sl    Sections along a map applied to higher-order Lagrangian
Mechanics. Noether's theorem}, 
Acta Applicandae Mathematicae {\bf 25} (1991), pp. 127--51. 

\medskip

[\Super] 
{\smc  Cari\~nena J.F. and Figueroa  H.}, {\sl Hamiltonian versus
Lagrangian formulation of supermechanics}, J. Phys. A: Math. Gen. 
{\bf 30} (1997), (to appear).

\medskip

[\Grapon] {\smc  Gr\`acia X. and  Pons J.M.}, {\sl On an evolution operator
connecting Lagrangian and Hamiltonian formalisms},
Lett. Math. Phys. {\bf 17} (1989), pp. 175--180.

\medskip

[\Marmuk] {\smc Marmo G. and  Mukunda N.}, {\sl Symmetries and constants
of the motion in Lagrangian mechanics: beyond point transformations},
Nuovo  Cim. {\bf A 92} (1986), pp. 1--12.

\medskip

[\MCS] {\smc  Mart\'{\i}nez E., Cari\~nena J.F.,   and  W. Sarlet},
{\sl Derivations of differential forms along the tangent bundle
projection},
 Diff. Geom. and Appl. {\bf 2} (1992), pp. 17--43.

\medskip

[\Pide] {\smc Pidello G. and Tulczyjew W.M.},  {\sl Derivations of
differential forms in jet bundles},  Ann. Math. Pura ed Aplicata {\bf 147}
(1987), pp. 249--265.

\medskip

[\Tulcz]  {\smc Tulczyjew W.M.},  {\sl The Lagrange Differential},
Bul. Acad. Pol. Sc., vol. {\bf XXIV} (12) (1976), pp. 1089--1096.

\medskip

[\Tulczdos]  {\smc Tulczyjew W.M.}, {\sl The Legendre transformation},
Ann. Inst. Henri Poincar\'e, vol. {\bf XXVII} (1) (1977), pp. 101--114.

}
\bigskip

\end

\bye